\documentclass[12pt]{article}
\usepackage{amssymb}
\usepackage{amsmath}

\usepackage{latexsym}
\usepackage{amsfonts}
\usepackage{times}

\newtheorem{theorem}{Theorem}[section]
\newtheorem{lemma}[theorem]{Lemma}
\newtheorem{corollary}[theorem]{Corollary}

\newenvironment{Proof}{\removelastskip\par\medskip
\noindent{\em Proof.} \rm}{\penalty-20\null\hfill$\square$\par\medbreak}

\begin{document}

\centerline{\textbf{\Large Generalized potentials}}

\centerline{\textbf{\Large on commutative hypergroups}}

\vspace{3mm} \centerline{by}

\vspace{3mm} \centerline{\textbf{Mubariz G. Hajibayov}}
%\footnote{partially supported by the INTAS
%grant Ref. Nr 06-1000015-5777.}}

\centerline{\it National Aviation Academy}
\centerline{\it and}
 \centerline{\it Institute of Mathematics and Mechanics, Baku, Azerbaijan}

\centerline{\it (hajibayovm@yahoo.com)}

\begin{abstract}
By the Hardy-Littlewood-Sobolev theorem the classical Riesz potential is bounded on
Lebesgue spaces. E. Nakai and H. Sumitomo \cite{NS} extended that theorem to the Orlicz
spaces. We introduce generalized potential operators on commutative hypergroups and
under some assumptions on the kernel we showed the boundedness of these operators from
Lebesgue space into certain Orlicz space. Our result is an analogue of Theorem 1.3 in \cite{NS}.
\end{abstract}

\vskip 5pt


\vskip 5pt

\section{Introduction}
For $0 < \alpha < n$, the operator
$$
R_\alpha f (x) =\int\limits_{R^n}|x-y|^{\alpha-n} f(y)dy
$$
is called a classical Riesz potential (fractional integral).

By the classical Hardy-Littlewood-Sobolev theorem, if $1 < p < \infty$ and $\alpha p < n$, then $R_\alpha f$
is a bounded operator from $L^p (R^n)$ into $L^q (R^n)$, where
$\dfrac {1}{q }= \dfrac {1}{p} - \dfrac {\alpha}{n}$
(see \cite{H}, \cite{Stein} ).

The Hardy-Littlewood-Sobolev theorem is an important result in the potential theory. There
are a lot of generalizations and analogues of that theorem. The boundedness of the Riesz potentials
on spaces of homogeneous type was studied in \cite{GV} and \cite{KK}. The Hardy-Littlewood-
Sobolev theorem was proved for the Riesz potentials associated to nondoubling measures in
\cite{KM}. In \cite{Gadj} and \cite{Haj}, generalized potential-type integral operators were considered and (p, q)
properties of these operators were proved. In \cite{N}, \cite{NS}, \cite{Nak}, \cite{HajSam} the Hardy-Littlewood-Sobolev
theorem was extended to Orlicz spaces for generalized fractional integrals. In \cite{G}, \cite{GM}, \cite{GO},  \cite{TX},
Riesz potentials on different hypergroups were defined and analogues of the Hardy-Littlewood-
Sobolev theorem were given for these operators.

In this paper, we define generalized fractional integrals on commutative hypergroups and
prove the analogue of Theorem 1.3 in \cite{NS} for the generalized fractional integrals on commutative
hypergroups. The obtained result is an extension of the Hardy-Littlewood-Sobolev theorem
given in \cite{G}, \cite{GM}, \cite{GO},  \cite{TX}, for Riesz potentials on different hypergroups

Let $K$ be a set. A function $\rho
:K\times K\rightarrow \left[ 0,\infty \right) $ is called
quasi-metric if:

\begin{enumerate}
\item
$\rho \left( x,y\right) =0\,\Leftrightarrow \;x=y ;$
\item
$\rho \left( x,y\right) =\rho \left( y,x\right) ;$
\item
there exists a constant $c\geq 1$ such that for every $x,y,z\in K$
$$
\rho \left( x,y\right) \leq c\left( \rho \left( x,z\right) +\rho
\left( z,y\right) \right) .
$$
\end{enumerate}
Let all balls $B(x,r)=\{y\in K: \, \rho(x,y)<r\}$ be $\lambda$-measurable and assume that the measure $\lambda$
fulfils the doubling condition
\begin{equation} \label{doubling}
0<\lambda B(x,2r) \le D\lambda B(x,r)<\infty.
\end{equation}
A space $(K,\rho,\lambda)$ which satisfies all conditions mentioned above is
called a space of homogeneous type (see \cite{CW}).
% We recall that the measure $\lambda$ on $K$ is called doubling if
%\begin{equation} \label{doubling}
%\lambda B(x,2r) \le D\lambda B(x,r)<\infty
%\end{equation}
%for  every ball $B(x,r)=\{y\in K: \, \rho(x,y)<r\}$, with $D>0$ not depending on $x$ and
%$r$.\\
%The set $K$ equipped with a pseudo-metric  $\rho$ and doubling measure $\lambda$ is
%called the space of homogeneous type and denoted by $(K,\rho,\lambda)$.

In the theory of locally compact groups there arise certain spaces which, though not groups, have some of the structure of groups. Often, the structure can be expressed in terms of an abstract convolution of measures on the space. %A hypergroup is a measure algebra which has many of the useful
% properties associated with the convolution measure algebra of a group.

A hypergroup  $(K,\ast)$ consists of a locally compact Hausdorff space $K$ together with a
bilinear, associative, weakly continuous convolution on the Banach space of all bounded
regular Borel measures on $K$ with the following properties:
\begin{itemize}
\item[1.] For all $x,y \in K$, the convolution of the point measures $\delta _x\ast \delta _y$ is a probability measure with compact support.

\item[2.] The mapping:  $(x,y)\mapsto supp(\delta _x\ast \delta _y)$ of $K\times K$ into $\mathcal{C}(K)$, is continuous where $\mathcal{C}(K)$ is the space of compact subsets of $K$ endowed with the   Michael topology, that is the topology generated by the subbasis of all
$$U_{V,W}=\{L\in\mathcal{C}(K): L\cap V\neq \varnothing, L\subset W \}$$
where $V,W$ are open subsets of $K$.

\item[3.] There exits an identity $e \in K$ such that
$\delta _e\ast \delta _x=\delta _x\ast \delta _e=\delta _x$ for all $x \in K$.

\item[4.] There exits a topological involution $\thicksim$ from $K$ onto $K$ such that $\left(x^\thicksim\right)^\thicksim=x$, for $x \in K$, with
$$
\left(\delta _x\ast \delta _y\right)^\thicksim =\delta _{y^\thicksim }\ast \delta _{x^\thicksim}
$$
and
$e \in supp(\delta _x\ast \delta _y)$ if and only if $ x=y^\thicksim$ for $x,y \in K$ where for any Borel set $B$, $\mu^\thicksim \left(B\right)=\mu\left(\{x^\thicksim :x \in B\}\right)$
(see \cite{J}, \cite{Sp}, \cite{BH}, \cite{L}).
\end{itemize}
If $\delta _x\ast \delta _y=\delta _y\ast \delta _x$  for all
$x, y\in K$, then the hypergroup $K$ is called commutative . It is known that every commutative hypergroup $K$ possesses a Haar measure which will be denoted by $\lambda$
(see \cite{Sp}). That is, for every Borel measurable function $f$ on $K$,
$$
\int \limits_K f(\delta _x\ast \delta _y)d\lambda (y)=\int \limits_K f(y)d\lambda (y) \,\,\, (x \in K).
$$
Define the generalized translation operators $T^x$, $x\in K$,
by
$$
T^xf(y) =\int \limits_Kfd(\delta _x\ast \delta _y)
$$
for all $y \in K$. If $K$ is a commutative hypergroup, then $T^xf(y)=T^yf(x)$ and
the convolution of two functions is defined by
$$
\left(f\ast g\right)(x)=\int \limits_KT^xf(y)g(y^\thicksim)d\lambda (y).
$$
Let $p>0$. By $L^{p} \left( K,\lambda \right) $  denote a class of all
$\lambda$-measurable functions $f:K\rightarrow \left( -\infty
,\,+\infty \right) $  with $\left\| f\right\| _{L^{p} \left( K,\lambda \right)} =\left(
\int\limits_{K}\left| f\left( x\right) \right| ^{p} d\lambda \left(
x\right) \right) ^{\frac{1}{p} }<\infty $.\\
A function $\Phi :[0,\infty]\rightarrow [0,\infty]$ is called an $N$-function if can be represented as
$$
\Phi \left(r\right)=\int \limits_{0}^{r}\phi \left(t\right)dt,
$$
where $\phi :[0,\infty]\rightarrow [0,\infty]$ is a left continuous nondecreasing function such that $\phi \left(0\right)=0$ and $\lim_{t\rightarrow \infty}\phi \left(t\right)=\infty.$\\
Let $\Phi$ is an $N$-function. Define the Orlicz space $L^\Phi \left(K,\lambda \right)$ to be the set of all locally integrable functions $f$ in $K$ for which
$$
\int \limits_K \Phi \left(\dfrac {|f\left(x \right)|}{\eta} \right)d\lambda\left(x \right)\ < \infty
$$
for some $\eta >0$. Here $L^\Phi \left(K,\lambda \right)$ is equipped with the norm
$$
\|f\|_\Phi=\inf \{\eta >0: \int \limits_K \Phi \left(\dfrac {|f\left(x \right)|}{\eta} \right)d\lambda\left(x \right)\leq 1\}.
$$
For $\Phi \left(r \right)=r^p, 1<p<\infty$, we have $L^\Phi \left(K,\lambda \right)=L^p \left(K,\lambda \right)$.\\
The notation $\chi _A(x)$ denotes the characteristic function of set $A$.\\
Define a function $\Lambda_x(y)=T^x\chi_{B(e,r)}(y^\thicksim)$.\\
We will assume that there exit constants $c_1>0$, $c_2>0$ and $c_3>0$ such that for every $x,y\in K$ and $r>0$
\begin{equation} \label{2}
supp \Lambda_x(\cdot)\subset B(x,c_1r)
\end{equation}
and
\begin{equation} \label{3}
\lambda B(x,r)T^x\chi_{B(e,r)}(y^\thicksim)\leq c_2\lambda B(e,r)\leq c_3r^N.
\end{equation}
As examples of hypergroups satisfying the conditions \eqref{2} and \eqref{3} can be taken Laguerre,Dunkl and Bessel hypergroups (see \cite{G}, \cite{GM}, \cite{GO}).\\
A non-negative function $a(r)$ defined on $[0,\infty)$  is called almost increasing
(almost decreasing), if there exist a constant $C>0$ such that
$$a(t_1)\leq C a(t_2)$$
for all $0< t_1<t_2<\infty$  ($0<t_2<t_1< \infty$, respectively).\\
For an increasing function $a:\left(0,\infty\right)\rightarrow\left(0,\infty\right),$ define
$$
I_a f(x)=\int \limits_KT^x\left(\dfrac{a(\rho(e,y))}{\rho(e,y)^{N}}\right)f(y^\thicksim)d\lambda(y)
$$
on the commutative hypergroup $(K,\ast)$ equipped with the quasi-metric $\rho$. If $a(r)=r^\alpha, 0<\alpha <N,$ then $I_a$ is the Riesz potential of order $\alpha.$\\
Now we formulate a main result of the paper.
 \begin{theorem}\label{theorem}
Let $(K,\ast)$ be a commutative hypergroup, with the quasi-metric $\rho$ and doubling Haar measure $\lambda$ satisfying the conditions \eqref{2} and \eqref{3}.
Assume that $1<p<\infty$ and $a=a(r)$ is non-negative almost  increasing function on $[0,\infty)$,
$\dfrac{a(r)}{r^\lambda}$ is almost decreasing for some $0<\lambda <\frac{N}{p}$
and
$$
 \int \limits_0^1\frac{a(t)}{t}dt <\infty.
$$
Then the operator $I_a$ is bounded from $L^p(K,\lambda)$ into the Orlicz space $L^\Phi(K,\lambda)$, where the $N$-function is defined by its inverse
$$
\Phi ^{-1} \left( r  \right) =\int\limits_{0}^{r}A\left( t^{-\frac{1}{N}}\right) t
^{-\frac{1}{p^\prime}} dt,
$$
where $
A\left( r\right) =\int\limits_{0}^{r}\frac{a \left( t\right)
}{t} dt$.
\end{theorem}
If we take $a(r)=r^\alpha, 0<\alpha <N,$ then we have Hardy-Littlewood-Sobolev theorem for the Riesz potential
$$
I_\alpha f(x)=\int \limits_KT^x\rho(e,y)^{\alpha -N}f(y^\thicksim)d\lambda(y)
$$
on the commutative hypergroup $(K,\ast)$.
\begin{corollary}
Let $(K,\ast)$ be a commutative hypergroup, with the quasi-metric $\rho$ and doubling Haar measure $\lambda$ satisfying the conditions \eqref{2} and \eqref{3}. If $0<\alpha <N,$ $1<p< \frac N\alpha$ and $\frac 1p-\frac 1q= \frac \alpha N,$ then $I_\alpha$ is a bounded operator from $L^p\left(K,\lambda\right)$ into $L^q\left(K,\lambda\right)$.
\end{corollary}
\section{Preliminaries}
Define  Hardy-Littlewood maximal function
$$
Mf(x)=\sup\limits_{r>0}\frac{1}{\lambda B(e,r)}\left(|f|\ast \chi _{B(e,r)}\right)(x)
$$
on commutative hypergroup $(K,\ast)$ equipped with the pseudo-metric $\rho$.
\begin{lemma}\label{max}
Let $(K,\ast)$ be a commutative hypergroup, with quasi-metric $\rho$ and doubling Haar measure $\lambda$. Assume that
%$\lambda B(x,r)>0$ for all $x \in K$, $r>0$ and
there exist constants $c_1>0$ and $c_2>0$ such that for every $x,y\in K$ and $r>0$
$$
supp \Lambda_x(\cdot)\subset B(x,c_1r)
$$
and
$$
\lambda B(x,r)T^x\chi_{B(e,r)}(y^\thicksim)\leq c_2\lambda B(e,r).
$$
Then
\begin{itemize}
\item[1)] The maximal operator $M$ satisfies a weak type $(1,1)$ inequality, that is, there exists a constant $C>0$ such that for every $f\in L^1(K,\lambda)$ and $\alpha>0$
$$
\lambda\{x:Mf(x)>\alpha\}\leq \frac{C}{\alpha}\int\limits_K|f(x)|d\lambda(x).
$$
\item[2)] The maximal operator $M$ is of strong type $(p,p)$, for $1<p\leq \infty$, that is,
\begin{equation} \label{maxbound}
\|Mf\|_{L^p(K,\lambda)}\leq C_p\|f\|_{L^p(K,\lambda)},
\end{equation}
for some constant $C_p$ and every $f\in L^p(K,\lambda)$.
\end{itemize}
\end{lemma}
\begin{Proof}
It is clear that there exists nonnegative integer\,  $m$   such that $c_1\leq 2^m$ and $\lambda B(x,c_1r)\leq D^m\lambda B(x,r)$, where $D$ is a constant on doubling condition \eqref{doubling}. Then we have
$$
Mf(x)=\sup\limits_{r>0}\frac{1}{\lambda B(e,r)}\int \limits_{K}T^x|f(y)|\chi _{B(e,r)}(y^\thicksim)d\lambda(y)
$$
$$
=\sup\limits_{r>0}\frac{1}{\lambda B(e,r)}\int \limits_{K}|f(y)|T^x\chi _{B(e,r)}(y^\thicksim)d\lambda(y)
$$
$$
\leq\sup\limits_{r>0}\frac{1}{\lambda B(e,r)}\int \limits_{B(x,c_1r)}|f(y)|T^x\chi _{B(e,r)}(y^\thicksim)d\lambda(y)
$$
$$
=\sup\limits_{r>0}\frac{1}{\lambda B(x,r)}\int \limits_{B(x,c_1r)}|f(y)|\frac{T^x\chi _{B(e,r)}(y^\thicksim)\lambda B(x,r)}{\lambda B(e,r)}d\lambda(y)
$$
$$
\leq c_2\sup\limits_{r>0}\frac{1}{\lambda B(x,r)}\int \limits_{B(x,c_1r)}|f(y)|d\lambda(y)\leq c_2D^m M_\rho f(x),
$$
where
$$
 M_\rho f(x)=\sup\limits_{r>0}\frac{1}{\lambda B(x,r)}\int \limits_{B(x,r)}|f(y)|d\lambda(y)
 $$
 is a maximal operator on $(K,\rho,\lambda)$.  It is well known that the maximal operator $M_\rho $ is of weak type $(1,1)$ and is bounded on $L^p(K,\lambda)$ (see \cite{CW}, \cite{ST}). This fact and the inequality $Mf(x)\leq c_2D^m M_\rho f(x)$ completes the proof.
 \end{Proof}
 
\section{Proof of Theorem \ref{theorem} }
We may suppose that $f(x)\geq 0$ and by the linearity of the operator $I_a$, it suffices to prove that $\|I_a
f\|_\Phi \leq C <\infty $ for $\|f\|_{L^p(K,\lambda)}\le 1$. accordance with Hedbergs trick,we split $I_{a } f\left( x\right)$ in the
standard way
$$
I_{a } f\left( x\right) =\int\limits_{B\left( e,r\right) }\frac{a \left( \rho\left(e,y\right)
\right) }{\rho\left(e,y\right) ^{N} } T^xf\left( y^\thicksim\right) d\lambda\left(y\right)
$$
$$
+\int\limits_{X
\setminus B\left( e,r\right) }\frac{a \left( \rho\left(e,y\right) \right)
}{\rho\left(e,y\right) ^{N} } T^xf\left( y^\thicksim\right) d\lambda (y) =\mathcal{A}_{r}(x) +\mathcal{B}_{r}(x).
$$
 Estimate $\mathcal{A}_{r}(x).$ Since $\frac{a(t)}{t^N}$ is almost decreasing,  we have
$$
\mathcal{A}_{r}(x)=\sum\limits_{k=0}^{\infty}\int\limits_{2^{-k-1}r\leq\rho\left(
e,y\right)<2^{-k}r}\frac{a \left( \rho\left(e,y\right) \right) }{\rho\left(e,y\right) ^{N}
} T^xf\left( y^\thicksim\right) d\lambda (y)
$$
$$
\leq C \sum \limits_{k=0}^{\infty}\frac{a\left(2^{-k-1}r \right) }{\left(2^{-k-1}r\right) ^{N} }
\int\limits_{2^{-k-1}r\leq\rho\left( e,y\right)<2^{-k}r} T^xf\left( y^\thicksim\right) d\lambda (y)\leq
C Mf\left(x\right)\sum \limits_{k=0}^{\infty}{a \left( 2^{-k-1}r \right) }
$$
$$
\leq C Mf\left(x\right)\sum \limits_{k=0}^{\infty}{\int \limits_{2^{-k-1}r}^{2^{-k}r} \dfrac
{a(t)}{t} dt} .$$ Therefore,
\begin{equation}\label{ar}
\mathcal{A}_{r}(x) \leq C A(r)Mf(x), \quad A(r)= \int \limits_0^r \frac{a(t)}{t}\,dt.
\end{equation}
Now estimate $\mathcal{B}_r(x).$ By the H\"older inequality  and the condition $\|f\|_{L^p(K,\lambda)}\le 1$,
  we obtain
$$
\mathcal{B}_r(x)\leq \left(\int \limits_{K\backslash B(e,r)}\left(T^xf(y^\thicksim)\right)^pd\lambda(y)\right)^{\frac{1}{p}}\left(\int \limits_{K\backslash B(e,r)}\left(\frac{a \left( \rho\left(e,y\right) \right) }{\rho\left(e,y\right) ^{N}
}\right)^{p'}d\lambda(y)\right)^{\frac{1}{p'}}
$$
$$
\leq \left(\int \limits_{K\backslash B(e,r)}\left(\frac{a \left( \rho\left(e,y\right) \right) }{\rho\left(e,y\right) ^{N}
}\right)^{p'}d\lambda(y)\right)^{\frac{1}{p'}}
$$
$$
=\left(\sum \limits_{k=0}^{\infty}\int \limits_{2^kr\leq \rho(e,y)<2^{k+1}r}\left(\frac{a \left( \rho\left(e,y\right) \right) }{\rho\left(e,y\right) ^{N}
}\right)^{p'}d\lambda(y)\right)^{\frac{1}{p'}}
$$
$$
\leq C \left(\sum \limits_{k=0}^{\infty}\left(\frac{a \left( 2^kr \right) }{\left(2^kr\right) ^{N}
}\right)^{p'}\int \limits_{\rho(e,y)<2^{k+1}r}d\lambda(y)\right)^{\frac{1}{p'}}
$$
$$
\leq C\left(\sum \limits_{k=0}^{\infty}\left(\frac{a \left( 2^kr \right) }{\left(2^kr\right) ^{N}
}\right)^{p'}\left(2^{k+1}r\right)^N\right)^{\frac{1}{p'}}
$$
$$
\leq C\left(\sum \limits_{k=0}^{\infty}\left(\frac{a \left( 2^kr \right) }{\left(2^kr\right) ^{\frac{N}{p}}
}\right)^{p'}\right)^{\frac{1}{p'}}
$$
$$
\leq C\left(\sum \limits_{k=0}^{\infty}\left( a \left( 2^kr \right)\right)^{p'}\int \limits_{2^kr}^{2^{k+1}r}\left(\frac{1 }{t ^{\frac{N}{p}}
}\right)^{p'}\frac 1tdt\right)^{\frac{1}{p'}}
$$
$$
\leq C\left(\sum \limits_{k=0}^{\infty}\int \limits_{2^kr}^{2^{k+1}r}\left(\frac{a(t) }{t ^{\frac{N}{p}}
}\right)^{p'}\frac 1tdt\right)^{\frac{1}{p'}}
$$
$$
= C\left(\int \limits_{r}^{\infty}\left(\frac{a(t) }{t ^{\frac{N}{p}}
}\right)^{p'}\frac 1tdt\right)^{\frac{1}{p'}}
$$
$$
\leq C\frac{a(r) }{r ^\beta}\left(\int \limits_{r}^{\infty}\left(t ^{\beta-\frac{N}{p}}
\right)^{p'}t^{-1}dt\right)^{\frac{1}{p'}}
$$
$$
\leq C\frac{a(r) }{r ^{\frac{N}{p}}}
$$
Therefore
\begin{equation}\label{br}
\mathcal{B}_r(x)\leq CA(r)r ^{-\frac{N}{p}}
\end{equation}
From \eqref{ar} and \eqref{br}, we have
$$
I_{a } f(x) \leq C \left( Mf\left( x\right) + r^{-\frac{N}{p } } \right) A\left(
r\right).
$$
Then
\begin{equation}\label{viaPhi}
I_{a } f(x)  \leq C\left[Mf(x)r^\frac{N}{p}+1\right]\Phi^{-1}\left( \frac{1}{r^N}\right)
\end{equation}
by Theorem 4.9 in \cite{HajSam}. If we choose  $r=[Mf(x)]^{-\frac{p}{N}}$, then the  inequality \eqref{viaPhi}
turns into
$$
I_{a } f(x)  \leq C \Phi^{-1}\left( [Mf(x)]^{p}\right)$$
and consequently ,$$\int \limits_{K}\Phi\left( \frac{I_a f(x)}{C}\right)d\lambda(x) \leq
\int \limits_{K}[Mf(x)]^{p}d\lambda(x)\leq 1,
$$
 where we have used \eqref{maxbound} and the fact that
 $\|f\|_{L^p(K,\lambda)}\leq 1$.
 Hence
$$
\|I_a f\|_\Phi\leq C,
$$
which completes the proof.

\textbf{Acknowledgement.} This work was supported by the Science Development Foundation
under the President of the Republic of Azerbaijan Grant EIF-2012-2(6)-39/10/1. The author
would like to express his thanks to Academician Akif Gadjiev for valuable remarks.

Mubariz G. Hajibayov\\
National Aviation Academy. Bine gesebesi,\\
 25-ci km, AZ1104, Baku, Azerbaijan\\
and\\
Institute of Mathematics and Mechanics of NAS of Azerbaijan,\\
9, B. Vahavzade str., AZ1141, Baku, Azerbaijan.

%E-mail: mubarizh@rambler.ru


\begin{thebibliography}{99}

%\bibitem{AS} C. Abdelkefi and M. Sifi, Dunkl translation and uncentered maximal operator
%on the Real Line, \textit{Int. J. Math. Math. Sci. 2007, Article ID 87808, 9 pp.
%doi:10.1155/2007/87808}


 \bibitem{BH} W. R. Bloom and H. Heyer, Harmonic analysis of probability measures on hypergroups, \textit{de Gruyter Stud. Math.,} vol. 20, Walter de Gruyter $\&$ Co., Berlin, 1995.

 %\bibitem{BX1} W. R. Bloom and Z. Xu, The Hardy-Littlewood maximal function for Ch$\acute{e}$bli-Trim$\grave{e}$che hypergroups, \textit{Applications of hypergroups and related measure algebras, }(Seattle, WA, 1993), \textit{Contemp. Math.,} \textbf{183}  (1995), 45-70.

%\bibitem{BX2} W. R. Bloom and Z. Xu, Maximal function on Ch$\acute{e}$bli-Trim$\grave{e}$che hypergroups, \textit{Infin. Dimens. Anal. Quantum Probab. Relat. Top.,} \textbf{3(3)} (2000),  403-434.


\bibitem{CW} R. R. Coifman and G. Weiss, Analyse harmonique non-commutative sur certains espaces homog\`{e}nes.(French) \textit{Lecture Notes in Math.,} \textbf{242}, Springer-Verlag, Berlin-New York, 1971

%\bibitem{CSc} W. C. Connett and A. L. Schwartz, A Hardy-Littlewood maximal ineguality for Jacobi type hypergroups, \textit{Proc. Amer. Math. Soc.,} \textbf{107} (1989), 137-143.

%\bibitem{CSt} J.l. Clerc and E.M. Stein, $L_p$-multipliers for noncompact symmetric spaces, \textit{Proc. Nat. Acad. Sci. U.S.A,} \textbf{71} (1974), 3911–3912.



   % \bibitem{D} C. F. Dunkl, The measure algebra of a locally compact hypergroup,
%\textit{Trans. AMS,} \textbf{179}  (1973), 331-348.


  \bibitem{Gadj} A. D. Gadjiev, On generalized potential-type integral operators, \textit{Functiones et Approximatio, UAM,} \textbf{25} (1997), 37-44.

   \bibitem{GV}A. E. Gatto and S. Vagi,  Fractional integrals on spaces of homogeneous type, \textit{Analysis and
Partial Differential Equations,}  (1990), 171-216.


 %\bibitem{GGHM} G. Gaudry, S. Giulini, A. Hulanicki, A. M. Mantero, Hardy-Littlewood maximal function on some solvable Lie groups, \textit{J. Austral. Math. Soc. Ser. A,} \textbf{45(1)} (1988), 78-82.



    \bibitem{G} V. S. Guliyev, On maximal function and fractional integral, associated with the Bessel differential operator, \textit{Math. Inequal. Appl.,} \textbf{6(2)} (2003), 317-330.


     %\bibitem{G2} V. S. Guliyev and M. Assal, On maximal function on the Laguerre hypergroup, \textit{Fract. Calc. Appl. Anal.,} \textbf{9(3)} (2006), 307-318.

     \bibitem{GM} V. S. Guliyev, Y. Y. Mammadov, On fractional maximal function and fractional integrals associated with the Dunkl operator on the real line, \textit{J. Math. Anal. Appl.,} \textbf{353} (2009), 449-459.

         \bibitem{GO} V. S. Guliyev, M.N. Omarova, On fractional maximal function and fractional integral
on the Laguerre hypergroup, \textit{J. Math. Anal. Appl.,} \textbf{340(2)} (2008), 1058-1068.

      \bibitem{Haj} M.G.Hajibayov $(L_p;L_q)$ properties of the potential-type integrals associated to non-doubling measures, \textit{Sarajevo J. Math.} \textbf{2 (15)} (2006), 173-180.

          %\bibitem{HajSamk}M. G. Hajibayov and S. G. Samko, Weighted estimates of generalized potentials
%in variable exponent Lebesgue spaces
%on homogeneous spaces \textit{ Operator Theory:
%Adv. and Appl.,} \textbf{210} (2010), 107-122.

       \bibitem{HajSam}  M. G. Hajibayov and S. G. Samko, Generalized potentials in variable exponent Lebesgue spaces
on homogeneous spaces, \textit{ Math. Nachr.,} \textbf{284(1)} (2011), 53-66.

  %\bibitem{HL} G. H. Hardy and J. E. Littlewood, A maximal theorem with function-theoretic applications,
%\textit{Acta Math.,} \textbf{54(1)} (1930), 81-116.

\bibitem{H} L. Hedberg, On certain convolution inequalities, \textit{Proc. Amer. Math. Soc.,} \textbf{36} (1972), 505-510.

    \bibitem{J}   R. L. Jewett,  Spaces with an abstract convolution of measures.  \textit{Adv. in Math.,} \textbf{18(1)} (1975),  1-101.

\bibitem{KK} V. M. Kokilashvili and A. Kufner, Fractional inteqrals on spaces of homogeneous type, \textit{Comment. Math. Univ. Carolinae,} \textbf{30(3)} (1989), 511-523.

    \bibitem{KM} V. Kokilashvili and A. Meskhi, Fractional integrals on measure spaces, \textit{Frac. Calc.
Appl. Anal.,} \textbf{4 (1)} (2001), 1-24.

\bibitem{L} M. Lashkarizadeh Bami, The semisimplicity of $L^1(K,w) $ of a weighted commutative
hypergroup $ K$, \textit{Acta Math. Sin. (Engl.ser),} \textbf{24(4)} (2008), 607-610.

\bibitem{N} E. Nakai, On generalized fractional integrals, \textit{Taiwanese J. Math.}  \textbf{5(3)} (2001), 587-602.


\bibitem{NS} E. Nakai and H. Sumitomo, On generalized Riesz potentials and spaces of some smooth
functions, \textit{Sci. Math. Japonicae,} \textbf{54 (3)} (2001), 463-472.

\bibitem{Nak} E. Nakai,  On generalized fractional integrals in the Orlicz spaces on spaces of homoge-
neous type,  \textit{Sci. Math. Japonicae,} \textbf{54(3)} (2001), 473 -487.




    % \bibitem{So} F. Soltani, Littlewood-Paley operators associated with the Dunkl operator on R, \textit{J. Funct. Anal.,} \textbf{221} (2005), 205-225.


         \bibitem{Sp} R. Spector, Measures invariantes sur les hypergroupes(French), \textit{Trans. Amer. Math. Soc.,} \textbf{239} (1978), 147-165.


\bibitem{Stein} E. Stein, Singular integrals and diferentiability properties of functions, \textit{Princeton Mathematical Series,} No. 30 Princeton University Press, Princeton, N.J. 1970.


 %\bibitem{Stem} K. Stempak, Almost everywhere summability of Laguerre series, \textit{Studia Math.,} \textbf{100(2)} (1991), 129-147.


\bibitem{ST} J. O. Str\"{o}mberg and A. Torchinsky, Weighted Hardy spaces. Lecture Notes in Math., \textbf{1381}, Springer-Verlag, Berlin, 1989.

 \bibitem{TX} S. Thangavelu and Y. Xu, Riesz transform and Riesz potentials for Dunkl transform,
  \textit{J. of Comput. and Appl. Math.,}  \textbf{199(1)} (2007), 181-195.


%\bibitem{TX2} S. Thangavelu and Y. Xu, Convolution operator and maximal function for Dunkl transform,
% \textit{arXive:math0412037v1,} (2004), 1-17, \textit{J. Anal. Math.,} \textbf{97} (2005), 25-55

%\bibitem{W} N. Wiener, The ergodic theorem, \textit{Duke Math. J.,} \textbf{5(1)} (1939),  1-18.


\end{thebibliography}
\end{document}